\documentclass[12pt,a4paper]{amsart}
\usepackage{amssymb, hyperref, enumerate, parskip}
\usepackage{latexsym}
\usepackage{amsthm}

\begin{document}

\theoremstyle{plain}

\newtheorem{thm}{Theorem}
\newtheorem{cor}[thm]{Corollary}
\newtheorem{lem}[thm]{Lemma}
\newtheorem{prop}[thm]{Proposition}
\newtheorem{ex}[thm]{Example}
\newtheorem{ax}[thm]{Axiom}
\newtheorem{ques}[thm]{Question}
\newtheorem{defn}{Definition}[section]
\newtheorem{rem}[thm]{Remark}
\newcommand{\qqed}{\hfill \ensuremath{\Box}}

\title[An interesting elliptic surface]{An interesting elliptic surface over an elliptic curve}

\author{Tetsuji Shioda}
\address{Department of Mathematics, Rikkyo University, Nishi-Ikebukuro, Toshima-ku, Tokyo 171, Japan}
\email{shioda@rikkyo.ac.jp}
\author{Matthias Sch\"{u}tt}
\address{Department of Mathematics, Harvard University, 1 Oxford Street, Cambridge, MA 02138, USA}
\email{mschuett@math.harvard.edu}
\dedicatory{To Michel Raynaud}
\date{\today}

\thanks{We thank the Dipartimento di Matematica "Frederico Enriques" of Milano University for kind hospitality and stimulating atmosphere. The first author was funded by Grant-in-Aid for Scientific Research No.~17540044. The second author gratefully acknowledges support from the MCRTN network "Arithmetic Algebraic Geometry" and from DFG under grant "Schu 2266/2-1". We thank the anonymous referee for helpful comments.}

\begin{abstract}
We study the elliptic modular surface attached to the commutator
subgroup of the modular group. This has an elliptic curve as base
and only one singular fibre. We employ an algebraic approach and
then consider some arithmetic questions. \end{abstract}

\maketitle

\section{Introduction}

Let $\Gamma'$ be the commutator subgroup of the modular group
$\Gamma=SL(2,\mathbb{Z})$. It is known that $\Gamma'$ is a
congruence subgroup of $\Gamma$ of index $12$ with
\[
-\begin{pmatrix}1 & 6\\0 &
1\end{pmatrix}\in\Gamma',\;\;-\begin{pmatrix}1 & 0\\0 &
1\end{pmatrix}\not\in\Gamma'.
\]
Let $S=S(\Gamma')$ denote the elliptic modular surface attached to
$\Gamma'$ in the sense of \cite{Sh-EMS}. This elliptic surface has
the remarkable property that it has only one singular fibre. In
Kodaira's notation \cite{Kodaira}, this fibre has type $I_6^*$
(see \cite[Ex.~5.9]{Sh-EMS}).

From an analytic viewpoint, $S$ has been studied by Stiller in
\cite{Stiller}. Here we follow a more algebraic approach. We then
consider some arithmetic questions. Throughout the paper, $k$ will denote  an algebraically closed fied of
characteristic $\neq 2,3$.

\section{The elliptic modular curve associated to $\Gamma'$}
\label{s:B}

Let $B=B(\Gamma')$ denote the modular curve attached to
$\Gamma'$. Since $\Gamma(6)\subset\Gamma'\cdot\{\pm \mathbf{1}\}$, this elliptic curve is
closely related to the modular curve $B(6)$ with level $6$-structure. In
fact,
\[
\Gamma' \cdot \{\pm \mathbf{1}\} = \Gamma^{2}\cap\Gamma^{3}
\]
with the subgroups of squares resp.~cubes in $\Gamma$
\cite{Newman}. Hence we first recall the classical formulae for
the elliptic curves with level $2$ and level $3$-structure (employing
Igusa's notation \cite{Igusa}). 

The elliptic curve with level $2$-structure is given in
\emph{Legendre form} with a parameter $\lambda\in k-\{0,1\}$:
\emph{Legendre form} with parameter $\lambda \neq 0,1$:
\begin{eqnarray*}
E_\lambda/k(\lambda):\;\;\; y^2 = x \,(x-1)\, (x-\lambda),\;\;\;\;\;\;\;\;\;
\end{eqnarray*}
\begin{eqnarray*}
j(E_\lambda) \; = \; 256\,
(\lambda^2-\lambda+1)^3/\lambda^2\,(\lambda-1)^2,
\end{eqnarray*}
\begin{eqnarray}\label{eqn:eta}
j(E_\lambda)-12^3 = \eta^2.
\end{eqnarray}

where $\eta =
8\,(\lambda+1)\,(\lambda-2)\,(2\lambda-1)/\lambda\,(\lambda-1).$

The elliptic curve with level $3$-structure takes the
\emph{Hessian form}
\[
E_\mu/k(\mu):\;\; X^3 + Y^3 + Z^3 - 3\,\mu\, XYZ = 0
\]
with parameter $\mu\in k,\, \mu^3\neq 1$. Transformation to
Weierstrass form leads to
\[
E_\mu/k(\mu): \;\; y^2 = x^3 - 27\,\mu\,(\mu^3+8)\,x + 54\,(\mu^6-20\mu^3-8).\\
\]
Setting $\xi = 3 \mu (\mu^3+8)/(\mu^3-1)$ gives
\begin{eqnarray}\label{eqn:xi}
j(E_\mu) = \xi^3.
\end{eqnarray}

Combining (\ref{eqn:eta}) and (\ref{eqn:xi}) in case
$j=j(E_\lambda)=j(E_\mu)$, we derive the defining equation of the modular curve $B$:
\begin{eqnarray}\label{eq:B}
B/k:\;\; \eta^2 = \xi^3 - 12^3.
\end{eqnarray}
We choose the point at $\infty$ as
origin of the group law and denote it by $o_B$.

\section{The elliptic modular surface attached to $\Gamma'$}

There is a remarkable elliptic surface over $B$ which has constant
discriminant and yet variable moduli (i.e. non-constant
$j$-invariant) with only one singular fibre. Namely, consider the elliptic curve
\begin{eqnarray}\label{eq:E}
E:\;\; y^2 = x^3 - 27\,\xi\,x - 54\,\eta
\end{eqnarray}
over the function field $k(B)$. It is immediate that $E$ has
discriminant
\[
\Delta = 2^6 3^9 (\xi^3 - \eta^2) = 6^{12}.
\]
Hence the associated elliptic surface $S$ has no singular fibre
over $B-\{o_B\}$. On the other hand, we have
\[
j(E) = \xi^3,
\]
so that $ord_{o_B}(j)=-6$. Since $ord_{o_B}(\Delta)=12$, the singular fibre over $o_B$ has type $I_6^*$ (cf.~\cite{Tate}).

As modular curve resp.~surface for $\Gamma'$, such $B$ and $S$ are naturally
unique. In particular, for $k=\mathbb{C}$, we have
\[
B\otimes\mathbb{C} = \overline{\mathbb{H}/\Gamma'},\;\;
S\otimes\mathbb{C} = S(\Gamma')
\]
with $\mathbb{H}$ the upper half plane. Note also that the
fundamental group of the pointed curve
$B(\mathbb{C})-\{o_B\}$ is isomorphic to $\Gamma'$.

However, $S$ is not the unique elliptic surface over $B$ with only
one singular fibre up to isomorphism. In fact, twisting its
defining equation (\ref{eq:E}) over the $2$-torsion points of $B$,
we obtain three more such surfaces. These are also modular and
mutually isomorphic (cf.~\cite{Stiller}). Among these surfaces,
$S$ is distinguished by the constant discriminant
$\Delta$. Any elliptic surface over an elliptic curve with constant discriminant and exactly one singular fibre is, up to isomorphism, obtained from $S$ via purely inseparable base change.

\section{The elliptic modular surface of level $6$}

We already pointed out that there is a close relation between
$\Gamma'$ and the principal congruence subgroup $\Gamma(6)$. In
this section, this will be made explicit.

In the notation of Section \ref{s:B}, consider the diagram of
extensions
$$
\begin{array}{ccccc}
&& \mathbb{Q}(\lambda,\mu)&&\\
& \nearrow && \nwarrow &\\
\mathbb{Q}(\lambda,\xi) &&&& \mathbb{Q}(\eta, \mu)\\
& \nwarrow && \nearrow &\\
&& \mathbb{Q}(\eta, \xi) &&
\end{array}
$$
It defines a commutative diagram of isogenies of elliptic curves
$$
\begin{array}{ccccc}
&& B(6) &&\\
& \swarrow && \searrow &\\
A &&&& \tilde{B}\\
& \searrow && \swarrow &\\
&& B &&
\end{array}
$$
where the $\swarrow$ are the duplication maps and the $\searrow$
are $3$-isogenies.

Then the base change of the elliptic curve $E$ over $k(B)$ from
$B$ to $B(6)$ gives, appropriately twisted, the elliptic modular
surface of level $6$. At the same time, we obtain the
identification of $\mathbb{Q}(B(6))=\mathbb{Q}(\lambda,\mu)$ with
$\mathbb{Q}(s,t)$ satisfying $s^2 = t^3 - 12^3$. This is perhaps
slightly more natural than the construction given in \cite{RS}.

\section{The cusp forms associated to $S$}

We recall some notable properties of elliptic modular surfaces
(cf.~\cite{Sh-EMS}). Let $\Gamma''\subset\Gamma$ of finite index
such that $-\mathbf{1}\not\in\Gamma''$. Consider the complex
elliptic modular surface $S=S(\Gamma'')$ attached to $\Gamma''$:

\begin{enumerate}[(i)]
\item The holomorphic $2$-forms on the elliptic modular surface
$S$ correspond to the cusp forms of weight 3 with respect to $\Gamma''$.
(This resembles the correspondence of holomorphic $1$-forms on the
modular curve $B(\Gamma'')$ with the cusp forms of weight 2 with respect
to $\Gamma''$.) In particular, the geometric genus $p_g$ of $S$ equals
the dimension of the $\mathbb{C}$-vector space of cusp forms
$\mathcal{S}_3(\Gamma'')$.

\item $S$ is extremal: The N\'eron-Severi
group has maximal rank $\rho(S)=h^{1,1}(S)$, while the
Mordell-Weil group has rank zero.
\end{enumerate}

For the commutator subgroup $\Gamma'$, we have $g=p_g=1$. Thus,
the Hodge diamond of $S=Y(\Gamma')$ reads
$$
\begin{array}{ccccc}
&&1&&\\& 1 && 1 &\\
1 && 12 && 1\\ & 1 && 1 &\\
&& 1 &&
\end{array}
$$

We shall now determine the precise cusp forms corresponding to our
model $S/\mathbb{Q}$, given by (\ref{eq:B}) and (\ref{eq:E}). For
some prime $\ell$, we consider the $\ell$-adic Galois
representations associated to $H_{et}^1(S, \mathbb{Q}_{\ell})$ and
to the transcendental lattice
\[
T_S=NS(S)^{\bot}\subset  H^2(S,\mathbb{Z}).
\]
By the above properties, both representations are two-dimensional
and correspond to some cusp forms with rational Fourier
coefficients $a_n$ resp.~$b_n$. Furthermore, $H_{et}^3(S,
\mathbb{Q}_{\ell})=H_{et}^1(S, \mathbb{Q}_{\ell})(1)$ by
Poincar\'e duality, so the corresponding $L$-series agree up to a
shift. The remaining cohomology is algebraic. It gives rise to
one-dimensional Galois representations which are easily
understood.

Recall Dedekind's eta-function
\[
\eta(\tau)=q \prod_{n\in\mathbb{N}}(1-q^n)\;\;\;(q=e^{2\pi
i\tau/24}).
\]
It is classically known that
\[
\mathcal{S}_2(\Gamma')=\mathbb{C}\, \eta(\tau)^4\;\;\text{ and
}\;\;\mathcal{S}_3(\Gamma')=\mathbb{C}\, \eta(\tau)^6.
\]
As a consequence, the cusp forms corresponding to $H_{et}^1(S,
\mathbb{Q}_{\ell})$ and $T_S$ are twists of the above forms,
depending on the chosen model $S/\mathbb{Q}$ (which implicitly
also includes the choice of the model $B/\mathbb{Q}$).

Note that $H_{et}^1(S, \mathbb{Q}_{\ell})=H_{et}^1(B,
\mathbb{Q}_{\ell})$. Hence we can use the known modularity of $B$
to deduce
\begin{eqnarray}\label{eq:wt2}
L(H_{et}^1(S, \mathbb{Q}_{\ell}), s) = L(\eta(\tau)^4, s).
\end{eqnarray}
Classically, this follows from the fact that $B$ has CM by $\mathbb{Q}(\sqrt{-3})$
(cf.~\cite[II.10.6]{Silverman}). Alternatively, since $B$ has
conductor $36$, it can be derived from the fact that
\[
\mathcal{S}_2(\Gamma_0(36))=\mathbb{C}\,\eta(6\tau)^4.
\]
In terms of $L$-series, there is no need to distinguish between these two cusp forms since $L(\eta(\tau)^4, s)=L(\eta(6\tau)^4, s)$. 

On the other hand, we consider $T_S$. Note that $S$ has good
reduction outside $\{2,3\}$. As this prescribes the ramification
of $H_{et}^{*}(S, \mathbb{Q}_{\ell})$, there are only a few
possible twists of $\eta(\tau)^6$. We have to take both
quadratic and quartic twists into consideration, since
$\eta(\tau)^6$ has complex
multiplication by $\mathbb{Q}(\sqrt{-1})$. The quartic twisting
can be achieved in terms of the corresponding Gr\"ossencharacter
of $\mathbb{Q}(\sqrt{-1})$ which has conductor $(2)$ and
$\infty$-type $2$ (cf.~\cite{S-CM}).

An observation of the twisting characters in question shows that
they are determined by their values at (the primes in
$\mathbb{Z}[\sqrt{-1}]$ above) $5$ and $13$. To determine the cusp
form corresponding to $T_S$, it thus suffices to know the Fourier
coefficients at $5$ and $13$. 

At a prime $p>3$, the coefficient
$b_p$ can be computed with the Lefschetz fixed point formula. We
use the modularity result of (\ref{eq:wt2}) and Poincar\'e duality as
mentioned, plus the fact that the one-dimensional
(i.e.~algebraic) representations involved are trivial. In other
words, $NS(S)$ is generated by divisors over $\mathbb{Q}$. In the present case, this follows from Tate's algorithm \cite{Tate}  which shows that all components of the $I_6^*$ fibre of $S$ are defined
over $\mathbb{Q}$. As a result,
the Lefschetz fixed point formula for $S$ reads
\[ \# S(\mathbb{F}_p) = 1 + 12p +
b_p + p^2 - (1+p) a_p.
\]
Since we know $a_p$, we can calculate $b_p$ from the number of
points of $S$ over $\mathbb{F}_p$. Counting points with a machine,
we obtain $b_5=-6$ and $b_{13}=10$. Up to the Euler factor at $3$,
this gives
\begin{eqnarray}
L(T_S, s) = L(\eta(\tau)^6, s).
\end{eqnarray}

\begin{thm}\label{thm:L}
Up to the Euler factors at $2$ and $3$, we have
\[
\zeta(S/\mathbb{Q}, s)= \frac{\zeta(s)\,\zeta(s-1)^{12}\,L(\eta(\tau)^6,  s)\,\zeta(s-2)}{L(\eta(\tau)^4, s)\, L(\eta(\tau)^4, s-1)}.
\]
\end{thm}

\section{The rank of $E$ depending on the characteristic}

We shall use Theorem \ref{thm:L} to determine the rank of
the elliptic curve $E$ of (\ref{eq:E}) in positive characteristic.
In other words, we are concerned with the Mordell-Weil rank and the
Picard number of the corresponding surface $S$ and look for
supersingular primes. Recall that $S$ has good reduction at the
primes $p>3$.

\begin{lem}\label{lem:rho}
Let $k$ be an algebraically closed field of characteristic $p>3$.
Then
\[
\rho(S/k)=\begin{cases}
12 & \text{ if } \;p\equiv 1\mod 4,\\
14 & \text{ if } \;p\equiv 3\mod 4.
                           \end{cases}
\]
\end{lem}
The leitmotif to prove the lemma is to consider the
$\zeta$-function of $S/\mathbb{F}_p$. This is obtained from
$\zeta(X/\mathbb{Q}, s)$ by considering the local Euler factors at
$p$. From the associated Gr\"ossencharacter, we know that the
factors corresponding to $\eta(\tau)^6$ have eigenvalues
\begin{eqnarray}\label{eq:pp}
\begin{cases} \pi^2, \bar\pi^2 & \text{ if } p\equiv 1\mod 4\;
\text{ splits as } p=\pi\bar\pi \,\text { in } \mathbb{Z}[2\sqrt{-1}],\\
p, -p & \text{ if } p\equiv 3\mod 4, \,p>3.
\end{cases}
\end{eqnarray}
As a consequence, Lemma \ref{lem:rho} in its entirety would follow
from the Tate Conjecture \cite{Tate-1}, but this is only known in
some few cases (cf.~\cite[Thm.~(5.6)]{Tate-Conj}). Nevertheless, the
first case of the lemma follows from (\ref{eq:pp}), since then the
eigenvalues are not $p$ times a root of unity, a necessity for
algebraic classes (cf.~\cite[App.~C, Ex.~10]{Sh-EMS}).

Now, by the relation of the Picard number and the Mordell-Weil rank (cf. \cite[Cor.~5.3]{Sh-MW}),
 Lemma \ref{lem:rho} is equivalent to the following:
\begin{eqnarray}\label{eq:MW}
\text{rank } E(k(B)) =\begin{cases}
0 & \text{ if } \;p\equiv 1\mod 4,\\
2 & \text{ if } \;p\equiv 3\mod 4, \,p>3.
                           \end{cases}
\end{eqnarray}
To prove the second case above, we use the fact that $S$ can be derived from an elliptic K3
surface via base change. The advantage of this approach is that
the Tate Conjecture is known for elliptic K3 surfaces.

Let $X$ be the elliptic K3 surface
over $\mathbb{P}^1$, given in Weierstrass form
\begin{eqnarray}\label{eq:XX}
E':\;\; y^2 = x^3 - 27\,(t^2+12^3)^3\, x - 54 \,t\, (t^2 + 12^3)^4.
\end{eqnarray}
This has singular fibres of type $I_2^*$ over $\infty$ and $IV^*$
over the two square roots of $-12^3$. The idea is to pull-back via
a base change of degree 3 which is ramified exactly above these
cusps. In terms of the affine modular curve $B$,
such a map is given by the following projection:
$$
\begin{array}{ccc} B & \to &
\mathbb{A}^1\\
(\xi, \eta) & \mapsto & \eta.
\end{array}
$$


The base change replaces $t$ in the Weierstrass equation (\ref{eq:XX}) of $X$ by $\eta$. Then we use the relation between $\eta$ and $\xi$ from the defining equation (\ref{eq:B}) of $B$. After minimalizing, we
obtain exactly the modular surface $S$ as in (\ref{eq:E}).  Projectively, the pull-back corresponds to the following commutative diagram:
\begin{eqnarray*}
S & \to & X\\
\downarrow && \downarrow\\
B & \to & \mathbb{P}^1.
\end{eqnarray*}

In other words, the elliptic curve  $E/k(B)$ is obtained from $E'/k(t)$ by this base change. Let us identify $E'=E$ below.
As a consequence, we have the injection $E(k(t)) \hookrightarrow
E(k(B))$ for any characteristic $\neq 2,3$. Hence, to
deduce the second claim of (\ref{eq:MW}), it suffices to prove the
corresponding statement for $E(k(t))$.

\begin{lem}\label{lem:MW}
If $k$ has characteristic
$p\equiv 3\mod 4,\, p>3$, then
\[
\text{rank } E(k(t))=2.
\]
\end{lem}

The proof starts by considering $X/\mathbb{Q}$ which is known to
be modular. This provides us with the $\zeta$-function of $X/\mathbb{F}_p$ for any $p>3$. Then we apply the (known) Tate
Conjecture.

Over $\mathbb{C}$, $X$ is a singular K3 surface. Since $E(\mathbb{C}(t))=0$ (cf.~Cor.~\ref{Cor:E(k(t))}), $X$ is extremal, and the discriminant of
$NS(X)$ is $-36$. By \cite[Ex.~1.6]{Livne}, the transcendental
lattice $T_X$ (as a Galois module) is associated to a newform of
weight $3$ with complex multiplication by $\mathbb{Q}(\sqrt{-1})$.
(By construction, this is exactly $\eta(\tau)^6$.) Hence, if
$p\equiv 3\mod 4,\, p>3$, the roots of the Euler factor at $p$ are
$p$ and $-p$ as in (\ref{eq:pp}).

We now consider the reduction of $X$ modulo $p$. Recall that the
Tate Conjecture is known for elliptic K3 surfaces in
characteristic $p>3$ \cite[Thm.~(5.6)]{Tate-Conj}. Due to the
eigenvalues, this predicts that the Picard number (over the
algebraic closure) increases by two upon reducing. Since the fibre
configuration stays unchanged, this implies Lemma \ref{lem:MW}. As
explained, Lemma \ref{lem:rho} follows.

In more detail, the above argument gives the following
\begin{cor}
Fix $p\equiv 3\mod 4,\, p>3$. Let $r\in\mathbb{N}$ and $q=p^r$.
Then
\[
\text{rank } E(\mathbb{F}_q(B))=\begin{cases} 1 & \text{ if $r$ is
odd,}\\ 2 & \text{ if $r$ is even.}\end{cases}
\]
In particular, the Tate Conjecture holds for $S/\mathbb{F}_q$.
\end{cor}

As the Tate conjecture implies the Artin-Tate conjecture (\cite{Milne}, cf.~\cite[App.~C, Ex.~10]{Sh-EMS}), we have also established the following:
\begin{prop} \label{prop:p^2}
For any $p\equiv 3\mod 4, \,p>3$, we have
$$
 \det NS(X) =\det NS(S) =-p^2.
$$
\end{prop}

\section{Mordell-Weil lattices.}
Let us recall two general formulas from the theory of Mordell-Weil lattices  \cite{Sh-MW}. For a moment, suppose $E$ is an elliptic curve over $K=k(B)$ ($B/k$ any curve) such that the associated elliptic surface $f:S \to B$ has at least one singular fibre. We identify the $K$-rational points $P \in E(K)$ with the sections $\sigma:B \to S$ and denote the image curve $Im (\sigma) \subset S$ by the symbol $(P)$.

{\bf (1) Height formula: } For any  $P \in E(K)$, we have 
\begin{eqnarray}\label{eq:height}
\langle P,P \rangle =2 \chi(S) + 2(PO) -\sum_v contr_{v}(P),
\end{eqnarray}
where $\chi(S)$ is the arithmetic genus of $S$ and $(PO)$ denotes the intersection number of the  section $(P)$ and the zero-section $(O)$ on the elliptic surface $S$ (\cite[Thm.~8.2]{Sh-MW}). The term
 $contr_{v}(P)$ is a local contribution at $v \in B$ such that the fibre at $v$ is a reducible singular fibre; its value is given by \cite[(8.16)]{Sh-MW}.

{\bf (2) Determinant formula: } Let $M$ be the Mordell-Weil lattice $E(K)/E(K)_{tor}$. By \cite[Thm.~8.7]{Sh-MW}, we have
\begin{eqnarray}\label{eq:det}
\det M/ |E(K)_{tor}|^2 = \pm \det NS(S)/ \det V_S
\end{eqnarray}
where $V_S$ denotes the trivial lattice generated by the zero-section $(O)$ and fibre components.

Going back to the previous situation, we have (regardless of the characteristic $\neq 2,3$)
\begin{eqnarray}\label{eq:V_S}
 V_S = <1> \oplus <-1> \oplus \, D_{10}[-1].
 \end{eqnarray}
Here the root lattice $D_{10}$  corresponds to the singular fibre of type
 $I_6^*$.  For any lattice $L$,  $L[n]$ denotes the lattice with
the pairing  multiplied by $n$. 
\begin{thm} \label{thm:MWL}
Let $E/K, K=k(B)$ be as in (4) with  $p=char (k) \neq 2,3$. Then 
\begin{enumerate}[(i)]
\item $E(K)=0$  if $  p=0 $  or $ p\equiv 1\mod 4.$ 
\item $E(K) \cong \mathbb{Z}^{\oplus 2}[\frac{p}{2}]$ 
if $p \equiv 3\mod 4$.
\end{enumerate}
\end{thm}
\begin{lem} \label{lem:free}
$E(k(B))$ is torsion-free.
\end{lem}
Indeed, if there is a torsion point $P \neq O$ in $ E(k(B))$, then  (\ref{eq:height}) implies the relation:
$$
0=\langle P,P \rangle =
2\chi+2(PO)-contr_{o_B}(P)
$$ 
where $\chi=1$ and 
$contr_{o_B}(P)=0,1$  or $ \frac{5}{2}$ for a singular fibre of type $I_6^*$. But this is impossible, 
since there is no integer $(PO) \geq 0$ satisfying the above relation.

Then the assertion $(i)$ of Theorem \ref{thm:MWL} follows in view of (\ref{eq:MW}).
To prove the assertion $(ii)$, note first that the rank two lattice $M=E(K)$ has $\det M= (p/2)^2$, by Prop.~\ref{prop:p^2}, (\ref{eq:det}) and (\ref{eq:V_S}).
The following lemma implies that $M$ is similar to the square lattice, which will complete the proof of $(ii)$.

\begin{lem} \label{lem:order4}
The elliptic surface  $S$ is endowed with an automorphism of order $4$:
\[
\phi: (x, y, \xi, \eta) \mapsto (-x, iy, \xi, -\eta),\;\;\;\;\; i=\sqrt{-1}. 
\]
It acts as an automorphism of order 4 on the trans\-cendental lattice $T_S$ in case $k=\mathbb{C}$,
and on $E(K)$ in the supersingular case $p \equiv 3\mod 4$.
\end{lem}
First consider the case $k=\mathbb{C}$. Let $\omega_S$ be the unique  
 holomorphic $2$-form on $S$ up to constants and suppose $\phi^* (\omega_S)=\alpha \omega_S$. 
Since the fibre involution $\phi^2$ acts as multiplication by $-1$,  we have $\alpha^2=-1$ and $\alpha=\pm i$. As $\omega_S$ is contained in $T_S \otimes \mathbb{C}$, this proves the first assertion.
 In this case, $H^2(S, \mathbb{Z})$ contains $NS(S) \oplus T_S$ as a finite index subgroup, and we have $NS(S)=V_S$ because $NS(S)/V_S \cong E(K)=0$ by  \cite[Thm.~1.3]{Sh-MW} and the above Theorem 6 $(i)$. The action of $\phi^*$ on $V_S$ is trivial except that it exchanges two simple components of the $I_6^*$ fibre. Hence it is of order 2.

Next  we pass to $S(p)$ which denotes the reduction of $S$ mod $p$.  
Let $k(p)=\bar{\mathbb{F}}_p$. In order to consider $S(p)/k(p)$, we work with \'etale cohomology, and choose $\ell \equiv 1\mod 4$ so that $i \in \mathbb{Q}_{\ell}$. 
Then $T_S\otimes\mathbb{Q}_{\ell}$ is identified as $\pm i$-eigenspace of $\phi^*$ on $H_{et}^2(S_{\bar{\mathbb{Q}}},\mathbb{Q}_{\ell})$   by what we have seen above.  Now assume $p \equiv 3\mod 4$. Then $NS(S(p))$ spans this cohomology group over $\mathbb{Q}_{\ell}$   by Lemma \ref{lem:rho} and the usual comparison theorem. Since the trivial lattice $V_S$ stays the same as in characteristic 0, $\phi^*$ acts on the Mordell-Weil group  $E(K) \cong NS(S(p))/V_S$ as an automorphism of order 4. Since the height pairing is invariant under automorphisms (cf.~\cite[Prop.~8.12]{Sh-MW}), $\phi^*$ is compatible with the lattice structure.  This completes the proof of Lemma \ref{lem:order4} and hence Theorem \ref{thm:MWL} $(ii)$.

\begin{cor}\label{Cor:E(k(t))}
We have 
\begin{enumerate}[(i)]
\item $E(k(t))=0$  if $  p=0 $  or $ p\equiv 1\mod 4,$ and 
\item $E(k(t)) \cong \mathbb{Z}^{\oplus 2}[\frac{p}{6}]$ 
if $p \equiv 3\mod 4$.
\end{enumerate}
\end{cor}
Indeed, $(i)$ is obvious since $E(k(t)) \subset E(k(B))$. As for $(ii)$, there is an injection of Mordell-Weil lattices $E(k(t)) \hookrightarrow E(k(B))[\frac{1}{3}]$ by \cite[Prop.~8.12]{Sh-MW}. Then equality follows from the determinant formula (\ref{eq:det}) applied to $X$.

\section{Similarity of lattices.}
We keep the notation $S, X$ to denote the elliptic surface over the elliptic curve $B$ and the elliptic K3 surface, defined by (4) and (9).
In this section, we  shall compare two kinds of lattices. On the one hand, we consider the transcendental lattices $T_X, T_S$. On the other hand, we consider 
the {\em supersingular reduction lattices} $L_X(p), L_S(p)$ which is defined for a supersingular prime $p$ as follows:
With respect to the natural injection 
$NS(X/\bar{\mathbb{Q}})\hookrightarrow NS(X/\bar{\mathbb{F}}_p)$, we take the orthogonal complement of its image:
$$
L'=NS(X/\bar{\mathbb{Q}})^\bot \subset NS(X/\bar{\mathbb{F}}_p).
$$
We define $L_X(p)=L'[-1]$,
and likewise for $S$. It was conjectured  in \cite{ShCR} that the positive-definite lattices $T_X$ and $L_X(p)$ are similar for a singular K3 surface $X$ and a supersingular prime $p$. This conjecture extends to any surface with transcendental lattice of rank two. 

In the present situation, it is easy to verify this:
\begin{prop}\label{prop:sim}
The lattices $T_S$, $T_X$, $L_S(p)$, $L_X(p)$ are all similar to the square lattice
 $L_0=\mathbb{Z}^{\oplus 2}$. In particular, the above-mentioned conjecture holds for $X$ and for $S$. More precisely, the scaling factor for each of these lattices
is determined as follows:
\begin{eqnarray*}
T_S=L_0[2], & \;& L_S(p)=L_0[2p],\\
T_X=L_0[6], &\;& L_X(p)=L_0[6p]. 
\end{eqnarray*}
\end{prop}
In fact, all the lattices $T_S, T_X, L_S(p), L_X(p)$ are rank two positive-definite lattices which admit a lattice automorphism of order 4 by Lemma \ref{lem:order4} and its analogue for $X$. Hence they are similar to the square lattice $L_0=\mathbb{Z}^{\oplus 2}$.

As for the scaling factor, it suffices to determine the discriminants of these lattices.  We have 
$$ \det T_X=  -\det NS(X) = -\det V_X =4\cdot 3^2= 6^2,$$
since $NS(X)/V_X \cong E(k(t))=0$, and similarly, $\det T_S=2^2$.
Next, $\det L_S(p)=(2p)^2$ and $\det L_X(p)=(6p)^2$ follow from Theorem \ref{thm:MWL} $(ii)$ resp.~Cor.~\ref{Cor:E(k(t))} $(ii)$ and the following 
\begin{lem}
Let $k=\bar{\mathbb{F}}_p$. 
The lattice $L_S(p)$ is isomorphic to the narrow Mordell-Weil lattice $E(k(B))^0$, and similarly, 
 $L_X(p)$ is isomorphic to $E(k(t))^0$.
\end{lem}
Indeed we have $NS(S)=V_S$ for $S/\bar{\mathbb{Q}}$, and the trivial lattice $V_S$ does not change under reduction mod $p$. Hence $L_S(p)$ agrees with the \emph{essential lattice} defined in \cite[\S8]{Sh-MW}. Thus the claim follows from \cite[Thm.~8.9]{Sh-MW}. The same argument applies to $X$.

\begin{rem}
In a recent preprint \cite{Shimada}, Shimada  studies transcendental lattices $T_X$ and supersingular reduction lattices $L_X(p)$ of singular K3 surfaces $X$. In the case at hand, i.e.~for the singular K3 surface $X$ in (\ref{eq:XX}) and $p\equiv 3\mod 4$ ($p>3$), Thm.~2 (L) in \cite{Shimada} gives an alternative proof that $L_X(p)=T_X[p]$.
\end{rem}

\section{Behaviour of $MW$ under base change}
\label{s:MW}

We start with a general question. Throughout we work over an algebraically closed field. Let $S$ be an elliptic surface over an elliptic curve $C$. Denote the Mordell-Weil rank of $S$ by $r$. Let 
\[
n_C: C \to C
\]
be multiplication by $n$. Write $S^{(n)}$ for the pullback of $S$ under $n_C$ and $r^{(n)}=$ rk $MW(S^{(n)})$.

\begin{ques}
How does $r^{(n)}$ behave as $n\to\infty$ ?
\end{ques}

By construction, $r^{(n)}/n^2$ is bounded. However, we will only be able to do better in the special situation where $S$ is extremal. Recall that $S$ is called extremal if its Picard number $\rho(S)$ is maximal and if its Mordell-Weil group $MW(S)$ is finite, i.e. $r=0$.

\begin{lem}\label{lem:extr}
If $S$ is extremal, then so is $S^{(n)}$.
\end{lem}

We will use the following invariants of an elliptic surface $S$ with elliptic base $C$: Let $\chi$ denote the arithmetic genus. In the case at hand, $\chi$ equals the geometric genus $p_g$. The Euler number is $e=12\,\chi$. This gives $b_2=12\,\chi+2$ and $h^{1,1}=10\,\chi+2$. Finally,
\[
\rho=r+2+\sum_{v\in C} (m_v-1)
\]
where $r$ is the Mordell-Weil rank as before and $m_v$ denotes the number of components of the fibre at $v\in C$.

Let us prove Lemma \ref{lem:extr} in characteristic zero. By assumption, $\rho(S)=h^{1,1}=10\,\chi+2$. Since $n_C$ is unramified, $e(S^{(n)})=n^2 e(S) = 12\, n^2 \chi$ and 
\begin{eqnarray*}
\rho(S^{(n)}) & = & r^{(n)} + 2 + n^2 \sum_{v\in C} (m_v-1)\\
& \geq & 2 + n^2 \sum_{v\in C} (m_v-1) = 2 + 10\, n^2 \chi.
\end{eqnarray*}

On the other hand, Lefschetz' bound reads 
\[
\rho(S^{(n)})\leq h^{1,1}(S^{(n)})=2 + 10\, n^2 \chi.
\]
Hence $\rho(S^{(n)})=h^{1,1}(S^{(n)})$ and $r^{(n)}=0$.

In case of positive characteristic, the assumption is $\rho(S)=b_2(X)$, and the same argument applies using Igusa's bound $\rho\leq b_2$.

\end{document}